# A Note on Noncommutative Even Square Rings

S K PANDEY[*]

**ABSTRACT**

In this note we find the least order of a noncommutative even square ring and we note that it is a nil ring having characteristic four. In order to prove the main result given in this note we mainly use suitable examples.

## Introduction

Recall that a ring $R$ is called an even square ring if $a^2 \in 2R, \forall a \in R$. We have studied some properties of even square rings in [1]. Each non-zero element of a commutative even square ring is a zero divisor provided it contains a unique non-zero element $a$ satisfying $a^2 = 2a = 0$.

It is worth to note that the zero square rings studied in [2] are just a particular class of even square rings. Each zero square ring is an even square ring however each even square ring is not a zero square ring. Clearly each zero square ring is a nil ring. But each even square ring is not a nil ring.

It has already been noted that each even square ring of characteristic two is a commutative ring and there does not exist any noncommutative even square ring of order four [1]. Therefore it is a natural question to ask that what is the least order of a noncommutative even square ring.

In this note we determine the least order of a noncommutative even square ring. We also note that the least order of a noncommutative even square ring of odd order is nine as it is seen in the case of non-even square rings.

The characteristic of a ring of order eight can be two, four or eight. But each even square ring of characteristic two is a commutative ring. Also if the order and characteristic of a ring are equal then the ring will be commutative. This leads to search whether or not there is an even square ring of characteristic four and order eight.

## RESULTS

**Proposition 1.** Let $R$ be an even square ring of characteristic four then $R$ is a commutative ring provided $ab + ba = 0, \forall a, b \in R$.

***Proof.*** Let $R$ be an even square ring of characteristic four. Clearly, $2a^2 = 0, \forall a \in R$. Let $a, b \in R$ then $(a+b)^2 = a^2 + b^2 + ab + ba$. This gives $2ab + 2ba = 0$. If $ab + ba = 0$ then $ab + ba = 2ab \Rightarrow ab = ba, \forall a, b \in R$

**Proposition 2.** Let $R$ be an even square ring of characteristic four such that $2ab = 2ba, \forall a, b \in R$ then $R$ is not necessarily a commutative even square ring.

***Proof.*** Please refer example1.

**Corollary.** Let $R$ be a ring such that $2ab + 2ba = 0, \forall a, b \in R$ then $R$ is not necessarily a commutative ring.

**Proposition 3.** The least order of a noncommutative even square ring of even order is eight and it is a nil ring having characteristic four.

***Proof.*** It has been already noted in [1] that there does not exist any noncommutative even square ring of order four.

The characteristic of a ring of order eight can be two, four or eight. But each even square ring of characteristic two is a commutative ring. Also if the order and characteristic of a ring are equal then the ring will be commutative. Thus there does not exist a noncommutative even square ring of order eight and characteristic two or eight.

Now in order to prove the result of this proposition we produce a noncommutative even square ring of order eight (please refer example 1) and characteristic four. This ring is a nil ring. It directly follows from the fact that each even square ring of order $2^n, n \in Z^+$ is a nil ring [1].

**Proposition 4.** The least order of a noncommutative even square ring of odd order is nine.

***Proof.*** Please refer example 2.

**Example 1.**

Let $R = \left\{ \begin{pmatrix} 0 & 0 \\ 0 & 0 \end{pmatrix}, \begin{pmatrix} 0 & 1 \\ 0 & 0 \end{pmatrix}, \begin{pmatrix} 0 & 2 \\ 0 & 0 \end{pmatrix}, \begin{pmatrix} 0 & 3 \\ 0 & 0 \end{pmatrix}, \begin{pmatrix} 2 & 0 \\ 0 & 0 \end{pmatrix}, \begin{pmatrix} 2 & 2 \\ 0 & 0 \end{pmatrix}, \begin{pmatrix} 2 & 1 \\ 0 & 0 \end{pmatrix}, \begin{pmatrix} 2 & 3 \\ 0 & 0 \end{pmatrix} \right\}$.

Clearly $R$ is a noncommutative even square ring of order eight under matrix addition and multiplication modulo four.

**Example 2.**

Let $R = \left\{ \begin{pmatrix} 0 & 0 \\ 0 & 0 \end{pmatrix}, \begin{pmatrix} 2 & 0 \\ 0 & 0 \end{pmatrix}, \begin{pmatrix} 0 & 2 \\ 0 & 0 \end{pmatrix}, \begin{pmatrix} 2 & 2 \\ 0 & 0 \end{pmatrix}, \begin{pmatrix} 4 & 0 \\ 0 & 0 \end{pmatrix}, \begin{pmatrix} 0 & 4 \\ 0 & 0 \end{pmatrix}, \begin{pmatrix} 4 & 4 \\ 0 & 0 \end{pmatrix}, \begin{pmatrix} 2 & 4 \\ 0 & 0 \end{pmatrix}, \begin{pmatrix} 4 & 2 \\ 0 & 0 \end{pmatrix} \right\}$

then $R$ is a noncommutative even square ring of order nine under matrix addition and multiplication modulo six. Hence the least order of a noncommutative even square ring of odd order is the same as it is seen in the case of noncommutative rings.

*Dept. of Mathematics,

Sardar Patel University of Police, Security and Criminal Justice,

Daijar-342304, Jodhpur, India.

E-mail: skpandey12@gmail.com